\documentclass[12pt,a4paper]{article}
\usepackage[latin1]{inputenc}
\usepackage[francais]{babel}
\usepackage{a4,amsmath,amssymb,amsfonts,amsthm}
\usepackage{amsrefs}
\usepackage[small,nohug,heads=vee]{diagrams}
\usepackage[all]{xy}

\begin{document}
\bibliographystyle{plain}


\def\mR{\M{R}}
\def\mZ{\M{Z}}
\def\mN{\M{N}}           
\def\mQ{\M{Q}}
\def\mC{\M{C}}
\def\mG{\M{G}}



\def\Spec{{\rm Spec}}
\def\rg{{\rm rg}}
\def\Hom{{\rm Hom}}
\def\Aut{{\rm Aut}}
 \def\Tr{{\rm Tr}}
 \def\Exp{{\rm Exp}}
 \def\Gal{{\rm Gal}}
 \def\End{{\rm End}}
 \def\det{{{\rm det}}}
 \def\Td{{\rm Td}}
 \def\ch{{\rm ch}}
 \def\che{{\rm ch}_{\rm eq}}
  \def\Spec{{\rm Spec}}
\def\Id{{\rm Id}}
\def\Zar{{\rm Zar}}
\def\Supp{{\rm Supp}}
\def\eq{{\rm eq}}
\def\Ann{{\rm Ann}}
\def\LT{{\rm LT}}
\def\Pic{{\rm Pic}}
\def\rg{{\rm rg}}
\def\et{{\rm et}}
\def\sep{{\rm sep}}
\def\ppcm{{\rm ppcm}}
\def\ord{{\rm ord}}
\def\Gr{{\rm Gr}}
\def\ker{{\rm ker}}
\def\rk{{\rm rk}}
\def\Lie{{\rm Lie}}


\def\beginProof{\par\underline{Preuve}: }
 \def\endProof{${\qed}$\par\smallskip}
 \def\pr{^{\prime}}
 \def\prpr{^{\prime\prime}}
 \def\mtr#1{\overline{#1}}
 \def\ra{\rightarrow}
 \def\mfp{{\mathfrak p}}

 \def\mQ{{\Bbb Q}}
 \def\mR{{\Bbb R}}
 \def\mZ{{\Bbb Z}}
 \def\mC{{\Bbb C}}
 \def\mN{{\Bbb N}}
 \def\mF{{\Bbb F}}
 \def\mA{{\Bbb A}}
  \def\mG{{\Bbb G}}
  \def\mP{{\Bbb P}}
 \def\CI{{\cal I}}
 \def\CE{{\cal E}}
 \def\CJ{{\cal J}}
 \def\CH{{\cal H}}
 \def\CO{{\cal O}}
 \def\CA{{\cal A}}
 \def\CB{{\cal B}}
 \def\CC{{\cal C}}
 \def\CD{{\cal D}}
 \def\CK{{\cal K}}
 \def\CL{{\cal L}}
 \def\CI{{\cal I}}
 \def\CM{{\cal M}}
  \def\CN{{\cal N}}
\def\CP{{\cal P}}
\def\CR{{\cal R}}
\def\CG{{\cal G}}
\def\CT{{\cal G}}
 \def\wt#1{\widetilde{#1}}
 \def\mod{{\rm mod\ }}
 \def\refeq#1{(\ref{#1})}
 \def\blb{{\big(}}
 \def\brb{{\big)}}
\def\mc{{{\mathfrak c}}}
\def\mcpr{{{\mathfrak c}'}}
\def\mcprpr{{{\mathfrak c}''}}
\def\ss{{\rm ss}}
\def\parf{{\rm parf}}
\def\P1{{{\bf P}^1}}
\def\cod{{\rm cod}}
\def\pr{\prime}
\def\prpr{\prime\prime}
\def\ss{\scriptstyle}
\def\OX{{ {\cal O}_X}}
\def\mpartial{{\mtr{\partial}}}
\def\inv{{\rm inv}}
\def\indlim{\underrightarrow{\lim}}
\def\prolim{\underleftarrow{\lim}}
\def\pprolim{'\prolim'}
\def\Pro{{\rm Pro}}
\def\Ind{{\rm Ind}}
\def\Ens{{\rm Ens}}
\def\without{\backslash}
\def\pbdb{{\Pro_b\ D^-_c}}
\def\qc{{\rm qc}}
\def\Com{{\rm Com}}
\def\an{{\rm an}}
\def\gfield{{\rm\bf k}}
\def\s{{\rm s}}
\def\dR{{\rm dR}}
\def\ari#1{\widehat{#1}}
\def\ul#1{\underline{#1}}
\def\sul#1{\underline{\scriptsize #1}}
\def\mou{{\mathfrak u}}
\def\ich{\mathfrak{ch}}
\def\cl{{\rm cl}}
\def\K{{\rm K}}
\def\R{{\rm R}}
\def\F{{\rm F}}
\def\L{{\rm L}}
\def\pgcd{{\rm pgcd}}
\def\rc{{\rm c}}
\def\N{{\rm N}}
\def\E{{\rm E}}
\def\H{{\rm H}}
\def\CHOW{{\rm CH}}
\def\A{{\rm A}}
\def\d{{\rm d}}
\def\Res{{\rm  Res}}
\def\GL{{\rm GL}}
\def\Alb{{\rm Alb}}
\def\alb{{\rm alb}}
\def\Hdg{{\rm Hdg}}
\def\Num{{\rm Num}}
\def\Irr{{\rm Irr}}
\def\Frac{{\rm Frac}}
\def\Sym{{\rm Sym}}
\def\indlim{\underrightarrow{\lim}}
\def\prolim{\underleftarrow{\lim}}
\def\red{{\rm red}}
\def\naive{{\rm naive}}
\def\ch{{\rm ch}}
\def\Td{{\rm Td}}
\def\T{{\rm T}}
\def\min{{\rm min}}
\def\slope{{\rm slope}}
\def\max{{\rm max}}
\def\min{{\rm min}}
\def\Sup{{\rm Sup}}
\def\Qb{\bar{\mQ}}
\def\mn{{\mu_n}}
\def\m2{{\mu_2}}
\def\TW{\{-1\}}
\def\BC{{\widetilde{\ch}}}
\def\ol#1{\overline{#1}}
\def\LocF{{\rm LocF}}
\def\LL{{\rm LL}}
\def\Coh{{\rm Coh}}
\def\Perf{{\rm Perf}}
\def\Sel{{\rm Sel}}
\def\twp{^{(p)}}
\def\Im{{\rm Im}}
\def\coLie{{\rm coLie}}


\def\RHom{{\rm RHom}}
\def\rRHom{{\mathcal RHom}}
\def\rHom{{\mathcal Hom}}
\def\dotimes{{\overline{\otimes}}}
\def\Ext{{\rm Ext}}
\def\rExt{{\mathcal Ext}}
\def\Tor{{\rm Tor}}
\def\rTor{{\mathcal Tor}}
\def\SP{{\mathfrak S}}
\def\perf{{\rm perf}}
\def\Bl{{\rm Bl}}
\def\ab{{\rm ab}}

\def\H{{\rm H}}
\def\D{{\rm D}}
\def\Del{{\mathfrak D}}
\def\Stab{{\rm Stab}}
\def\Div{{\rm Div}}
\def\Ver{{\rm Ver}}
\def\insep{{\rm insep}}

 \newtheorem{theor}{Th\'eor\`eme}[section]
 \newtheorem{prop}[theor]{Proposition}
 \newtheorem{propdef}[theor]{Proposition-D\'efinition}
 \newtheorem{cor}[theor]{Corollaire}
 \newtheorem{cordef}[theor]{Corollaire-Définition}
 \newtheorem{lemme}[theor]{Lemme}
 \newtheorem{lemma}[theor]{Lemma}
 \newtheorem{sublem}[theor]{sub-lemma}
 \newtheorem{defin}[theor]{D\'efinition}
 \newtheorem{conj}[theor]{Conjecture}
 \newtheorem{rem}[theor]{Remarque}
 \newtheorem{comp}[theor]{Complément}

 \parindent=0pt
 \parskip=5pt

\author{D. R\"ossler}

\title{Le groupe de Selmer des isogénies de hauteur un}

\maketitle

\begin{abstract}On montre que le groupe de Selmer d'une isogénie de hauteur un entre deux  variétés abéliennes définies sur le corps de fonctions d'une variété quasi-projective et lisse $V$ sur un corps parfait $k_0$ de caractéristique $p>0$ peut être plongé dans le groupe des homomorphismes entre deux fibrés vectoriels naturels sur $V$. 
\end{abstract}

\section{Introduction}

Soit $V$ un schéma intègre, qui est quasi-projectif et lisse sur un corps parfait $k_0$ de caractéristique $p>0$. Soit $K$ le corps de fonctions de $V$. 
Soit 
$\pi:\CA\to V$ (resp. $\rho:\CB\to V$) un schéma semiabélien sur $V$.  On notera 
$A$ (resp. $B$) sa fibre au-dessus du point générique de $V$. Soit $V_\ab\subseteq V$ l'ensemble 
des points $v\in V$ tels que $\CA_v$ est une variété abélienne sur $\kappa(v)$. On suppose que $V_\ab$ est un ouvert non vide et que son complément $E:=V\backslash V_\ab$ (muni de sa structure de sous-schéma fermé réduit) est un diviseur à croisements normaux au-dessus de $k_0$, au sens de \cite{Illusie-Reduction}. 

Si $S$ est un schéma de caractéristique $p$, on notera $$F_S:S\to S$$ l'endomorphisme de Frobenius absolu de $S$. Si $G$ est un schéma en groupes sur un schéma $S$, on 
écrira $\epsilon_{G/S}:S\to G$ pour la section nulle et $\omega_G=\omega_{G/S}:=\epsilon_{G/S}^*(\Omega_{G/S})=\Lie(G)^\vee$ pour l'algèbre de coLie de $G$ sur $S$. 

Soit $\iota:\CA\to\CB$ une isogénie de hauteur un. 
On rappelle que $\iota$ est par définition de hauteur un si $\ker\,\iota$ est un schéma en groupes de hauteur un (par définition, cela veut dire que $\ker\,\iota$ est fini et plat et que $F_{\ker\,\iota}=0$).  On écrira $\Gamma:=\ker\,\iota.$

En préparation de la formulation de notre résultat principal, on rappelle la proposition suivante, dont une variante est démontrée dans \cite[§2, "The second exact sequence"]{AM-Duality} (voir aussi \cite[III.3.5.6]{Milne-Arith}). 
Nous en donnerons une démonstration dans la sous-section \ref{ssMP}. 
On notera $K^{[1]}$ le corps $K$ vu comme une $k_0$-algèbre via l'application naturelle de 
$k_0$ dans $K$ précomposée par $F_{k_0}$.

\begin{prop}[Artin-Milne]
Il existe un homomorphisme de groupes canonique
$$
\Phi_{\Gamma_K}=\Phi_{\iota_K}:H^1(K,\Gamma_K)\hookrightarrow\Hom_K(F_{K}^*(\omega_{\Gamma_K}),\Omega_{K^{[1]}/k_0})
$$
qui est injectif.\label{MilneProp}
\end{prop}
Ici $H^1(K,\Gamma_K)$ est l'ensemble des classes d'isomorphisme de torseurs sous $\Gamma_K$ et sur 
$\Spec\,K$. On rappelle que cet ensemble a une structure canonique de groupe commutatif. 

Voici une description explicite de l'application $\Phi_{\Gamma_K}$. 

On remarque d'abord 
qu'on a isomorphisme canonique $\Omega_{K/k_0}\simeq\Omega_{F_K}$. Celui-ci provient de la suite exacte canonique 
$$
F_K^*(\Omega_{K/k_0})\to\Omega_{K^{[1]}/k_0}\to\Omega_{F_K}\to 0, 
$$
où la première flèche s'annule. 

Soit $\lambda:P\to\Spec\,K$ un torseur sous $\Gamma_K$. Le torseur 
$F_K^*(P)$ sous $F^*_K(\Gamma_K)$ est trivial d'après \cite[III.3.5.7]{Milne-Arith} et cette trivialisation 
est unique puisque $F^*_K(\Gamma_K)$ est un schéma en groupes infinitésimal. 
Il existe ainsi une unique flèche $\sigma:\Spec\,K\to P$ telle que 
$\lambda\circ\sigma=F_K$. La flèche $\sigma$ induit un morphisme 
de différentielles
$$
\sigma^*(\omega_{P/K})\to 
\Omega_{F_K}
$$
qui est par définition $\Phi_{\Gamma_K}(P)$, via les identifications $\Omega_{K^{[1]}/k_0}\simeq\Omega_{F_K}$ et $\sigma^*(\omega_{P/K})\simeq F_K^*(\omega_{\Gamma_K})$. 

Dans \cite[Remark (2.7)]{AM-Duality}, on trouve une description 
légérement différente de l'application $\Phi_{\Gamma_K}$. Les deux descriptions coïncident mais nous ne démontrerons pas cette coïncidence. 

On rappelle aussi l'existence de la suite "de Kummer"
\begin{equation}
0\to\Gamma_K\to A(K)\stackrel{\iota_K\,\,\,}{\to} B(K)\stackrel{\delta_{\iota_K}\,\,\,\,\,}{\to} H^1(K,\Gamma_K)\stackrel{t}{\to} H^1(K,A)
\label{KSeq}
\end{equation}
associée à $\iota_K$ (cf. \cite[Appendix C.4]{HS-Dio}). On a une suite semblable pour le changement 
de base de $A$, $B$ et $\iota_K$ à n'importe quel $K$-schéma $S$. 

On rappelle aussi la définition suivante (cf. \cite[Appendix C.4.1]{HS-Dio}):

\begin{defin}
Le groupe de Selmer $$\Sel(K,\iota_K)=\Sel_V(K,\iota_K)\subseteq H^1(K,\Gamma_K)$$ relativement à $V$ est 
l'ensemble des classes d'isomorphismes de torseurs $\lambda:P\to\Spec\,K$ sous $\Gamma_K$ tels que on 
a $P_{K_v}\in\delta_{K_v}(B(K_v))$ pour tout point $v\in V$ de codimension un. 
\end{defin}
On dit qu'un point $v\in V$ est de codimension un si sa clôture de Zariski 
est un sous-schéma fermé de codimension un. Un point 
de codimension un dans $V$ définit une valuation discrète dans $K$ et on note $K_v$ la complétion de $K$ le long de cette valuation.

Notons  $V^{[1]}$ le schéma $V$ vu comme un $k_0$-schéma via la flèche de 
$V$ dans $\Spec\, k_0$ qui est la composition de la flèche structurale avec  $F_{k_0}$. 
On écrira $\Omega_{V^{[1]}/k_0}(\log\, E)$ pour le faisceau des formes différentielles de $V^{[1]}\backslash E$ sur $k_0$ à singularités logarithmiques le long de $E$.  Le faisceau 
$\Omega_{V^{[1]}/k_0}(\log\, E)$ est cohérent et localement libre. Voir \cite[Intro.]{Illusie-Reduction} pour la définition et plus de détails 
sur cette notion. 

On remarque qu'on a une flèche naturelle 'restriction à la fibre générique' 
\begin{equation}
\rho:\Hom_V(F_{V}^*(\omega_\Gamma),\Omega_{V^{[1]}/k_0}(\log\,E))\to\Hom(F_K^*(\omega_{\Gamma_K}),\Omega_{K^{[1]}/k_0})
\label{IM}
\end{equation}
qui est injective puisque $V$ est intègre et les faisceaux cohérents en jeu sont sans torsion. 

Nous sommes finalement en mesure d'énoncer le résultat dont la démonstration est l'objet de la présente note. 

\begin{theor}
Soit $m\geq 3$. 
Supposons que $(p,m)=1$ et que $A[m](K)\simeq(\mZ/m\mZ)^{2g}$. 
Alors $$\Phi_{\iota_K}(\Sel_V(K,\iota_K))\subseteq\Im(\rho).$$ \label{THprinc}
\end{theor}

Au vu du Théorème \ref{THprinc} et de l'injectivité de $\rho$, on a donc

\begin{cor}
Soit $m\geq 3$. On suppose que $(p,m)=1$ et que $A[m](K)\simeq(\mZ/m\mZ)^{2g}$. 
Il existe un homomorphisme canonique
$$
\phi_{\iota_K}:\Sel(K,\iota_K)\hookrightarrow\Hom_V(F_{V}^*(\omega_\Gamma),\Omega_{V^{[1]}/k_0}(\log\,E))
$$
qui est injectif.
\label{corha}
\end{cor}

On peut en déduire le 

\begin{cor}[du Corollaire \ref{corha}]
Supposons que $\dim(V)=1$. Il existe alors un homomorphisme canonique
$$
\phi_{\iota_K}:\Sel(K,\iota_K)\hookrightarrow\Hom_V(F_{V}^*(\omega_\Gamma),\Omega_{V^{[1]}/k_0}(E))
$$
qui est injectif.
\label{corcourbes}
\end{cor}
Ici $\Omega_{V^{[1]}/k_0}(E):=\Omega_{V^{[1]}/k_0}\otimes\CO(E)$, où $\CO(E)$ est le fibré en droites 
associé à $E$. La démonstration du fait que le Corollaire \ref{corha} implique le Corollaire \ref{corcourbes} est une variante de la fin de la démonstration du Théorème \ref{THprinc} (voir après le Lemme \ref{lemetext}) et nous la laissons en exercice au lecteur.

\begin{rem}\rm Il est probable que l'hypothèse $A[m](K)\simeq(\mZ/n\mZ)^{2g}$ est inutile. 
Cette hypothèse est ici forcée par la nature de la démonstration, qui se fonde 
sur l'existence d'espaces de modules pour les variétés abéliennes munies de certaines structures de niveau. Elle ne joue pas un rôle important dans les applications. \end{rem}

On se demandera à raison s'il est toujours possible d'associer à une variété abélienne 
$A/K$ des modèles $V$ et $\CB$ comme plus haut, 
tels que $V$ soit une variété {\it projective} sur $k_0$. Nous répondons à cette question 
dans la remarque \ref{goodhyp} plus bas. 

La démonstration du Théorème \ref{THprinc} se fait en deux étapes. 

On démontre d'abord (cf. sous-section \ref{ssFro}) le cas particulier correspondant au morphisme de Frobenius 
relatif et au sous-groupe du groupe de Selmer donné par les torseurs provenant de points $K$-rationnels sur $B$. Dans ce cas ci, et avec la restriction supplémentaire que $\dim(V)=1$, un énoncé du même type se trouve déjà dans \cite{Rossler-On-the-group}, dont on 
reprendra la méthode.  Pour se ramener à (la généralisation de) l'énoncé démontré dans \cite{Rossler-On-the-group}, on doit démontrer la compatibilité entre 
deux flèches définies de manière différentes.  C'est ce qui est fait dans le Lemme \ref{alphalem} et le Corollaire \ref{corai}.

La deuxième étape (cf. sous-section \ref{ssFD}) du théorème consiste à se ramener 
au cas particulier démontré dans la première étape. Nous menons à bien cette réduction en remarquant que le morphisme $\Phi_{\iota_K}$ est compatible de manière naturelle avec la composition d'isogénies de hauteur un et en remarquant que 
si un torseur se trouve dans le groupe de Selmer alors il est donné par un point rationnel si on accepte de passer à une extension finie séparable d'un certain type de $K$. Ceci est une conséquence du théorème d'approximation de Greenberg.

\section{Démonstration du Théorème \ref{THprinc}}

\subsection{Démonstration de la Proposition \ref{MilneProp}}

\label{ssMP}

Il nous faut montrer que $\Phi_{\Gamma_K}$ est un homomorphisme de groupes et que cet homomorphisme est injectif. 

Nous commençons pas démontrer la première assertion, à savoir que $\Phi_{\Gamma_K}$ un homomorphisme de groupes. Soit 
$$\lambda_{P}:P\to\Spec\,K$$ et $$\lambda_{P'}:P'\to\Spec\,K$$ deux torseurs sous 
$\Gamma_K$. On considère $T:=P\times_{\Spec\,K} P'$ et on écrit $\lambda_T:T\to\Spec\,K$ pour le morphisme structural. Le schéma $T$ est un torseur sous 
$\Gamma_K\times_K\Gamma_K$. On sait (\cite[III.3.5.7]{Milne-Arith}) qu'il existe des uniques morphismes 
$\sigma_P:\Spec\, K\to P$ et $\sigma_{P'}:\Spec\,K\to P'$ tels que $$\lambda_{P}\circ\sigma_P=\lambda_{P'}\circ\sigma_{P'}=F_K.$$ 
Comme $\Gamma_K\times_K\Gamma_K$ est aussi un schéma en groupes infinitésimal, on 
sait aussi qu'il existe une unique flèche $\sigma_T:\Spec\,K\to T$ telle que 
$\lambda_T\circ\sigma_T=F_K$ et ainsi on doit avoir
$$
\sigma_T=\sigma_P\times_K\sigma_{P'}.
$$
La somme $\lambda_Q:Q\to\Spec\,K$ de $P$ et $P'$ dans $H^1(K,\Gamma_K)$ est par définition 
le quotient de $T$ par l'action de $\Gamma_K$ sur $T$ donnée en coordonnées par 
la formule 
$$
\gamma((p,p')):=(\gamma(p),\gamma^{-1}(p'))
$$
(pour la représentabilité du quotient, cf. par ex. \cite{Raynaud-Passage}).  
Soit $q:T\to Q$ le morphisme quotient. 
On sait à nouveau qu'il existe une unique flèche $\sigma_Q:\Spec\,K\to Q$ telle que 
$\lambda_Q\circ\sigma_Q=F_K$ et on voit donc que $q\circ\sigma_T=\sigma_Q.$ 
On a ainsi un diagramme commutatif

\begin{equation}
\xymatrix{
\Spec\,K\ar@/^3pc/[rrrr]^{=}\ar[r]\ar[ddr]_{\sigma_P} \ar@/_5pc/[dddrr]_{F_K}& F_K^*(\Gamma_K)\times_K F_K^*(\Gamma_K)\ar[d]^{\sim}\ar[rr]^{\,\,\,\,\,\,(\gamma,\gamma')\mapsto \gamma\cdot 
\gamma'} & & F_K^*(\Gamma_K)\ar[d]^{\sim}&\Spec\,K\ar[l]\ar[ddl]^{\sigma_Q}\ar@/^5pc/[dddll]^{F_K}\\
&F_K^*(T)\ar[rr]^{F_K^*(q)}\ar[d] & & F_K^*(Q)\ar[d]&\\
&T\ar[rr]^{q}\ar[dr] & & Q\ar[dl]&\\
& &\Spec\,K & &
}
\label{ddone}
\end{equation}

La première ligne du diagramme \refeq{ddone} donne une application 
d'algèbres de Lie
$$
F_K^*(\Lie(\Gamma_K))\oplus F_K^*(\Lie(\Gamma_K))\stackrel{(x,x')\mapsto x+x'}{\to}F_K^*(\Lie(\Gamma_K))
$$
dont le dual est l'application
$$
F_K^*(\omega_{\Gamma_K})\stackrel{x\mapsto (x,x')}{\to}
F_K^*(\omega_{\Gamma_K})\oplus F_K^*(\omega_{\Gamma_K}).
$$
On en déduit la première assertion.

Pour la deuxième assertion, il nous faut montrer que si 
$$\lambda_P:P\to\Spec\,K$$
est un $\Gamma_K$-torseur et que 
$\Phi_{\Gamma_K}(P)=0$ alors $P$ est un $\Gamma_K$-torseur trivial. 
Soit $s$ le point fermé de $P$ qui est l'image de $\sigma_P$, vu 
comme sous-schéma fermé réduit. 
Notons $\iota:s\to P$ le morphisme d'immersion et $\sigma_s:\Spec\,K\to s$
le morphisme tel que $\iota\circ\sigma_s=\sigma_P$. 

Le morphisme $\sigma_s$ et le morphisme structural de $s$ donnent des extensions de corps $$K\,|\,\kappa(s)\,|\,K^p=K.$$ Puisque le $\Gamma_K$-torseur $P$ est topologiquement réduit au point $s$, il est trivial ssi l'extension $\kappa(s)|K^p$ est triviale, autrement dit si $\kappa(s)=K^p$. 
Rappelons maintenant que $\Phi_{\Gamma_K}(P)$ est provient de l'application naturelle
$\sigma_P^*(\Omega_{P/K})\to\Omega_{F_K}$. Cette application se factorise
de la manière suivante:
$$
\sigma_P^*(\Omega_{P/K})=\sigma_s^*(\iota^*(\Omega_{P/K}))\to\sigma_s^*(\Omega_{s/K})\to\Omega_{F_K}
$$
où la première flèche est l'image réciproque par $\sigma_s$ de la flèche 
surjective 
$$
\iota^*(\Omega_{P/K}))\to\Omega_{s/K}
$$
induite par $\iota$. On voit donc que l'application $\Phi_{\Gamma_K}(P)$ est nulle 
ssi l'application $$\sigma_s^*(\Omega_{s/K})\to\Omega_{F_K}$$ est nulle. 
Remarquons maintenant qu'on a une suite exacte canonique
$$
\sigma_s^*(\Omega_{s/K})\to\Omega_{F_K}\to\Omega_{\sigma_s}\to 0
$$
On voit ainsi que $\Phi_{\Gamma_K}(P)$ est nulle ssi 
$\rk(\Omega_{F_K})=\rk(\Omega_{\sigma_s})$. Selon \cite[Th. 26.5, p. 202]{Matsumura-Commutative} 
on a 
$$
p^{\rk(\Omega_{\sigma_s})}=[K:\kappa(s)]
$$
et aussi
$$
p^{\rk(\Omega_{F_K})}=[K:K^p]
$$
ce qui implique pour finir que $\Phi_{\Gamma_K}(P)$ est nulle ssi 
$\kappa(s)=K^p$, ce qu'on voulait démontrer.

\begin{comp}
Supposons donné un diagramme commutatif
\begin{center}
\hskip2.4cm
\xymatrix{
\CA\ar[r]^{\iota}\ar[d]^{a} & \CB\ar[d]^{b}\\
\CA_1\ar[r]^{\iota_1} &\CB_1
}
\end{center}
où $\CA_1$ et $\CB_1$ sont des schémas semiabéliens sur $V$,
$\iota_1$ est une isogénie de hauteur un et $a$ et $b$ sont des morphismes de schémas en groupes. Alors 
on a $a(\ker\,\iota)\subseteq\ker\,\iota_1$ et pour tout $w\in F_K^*(\omega_{\ker\,\iota_{1,K}})$ et 
$P\in H^1(K,\ker\,\iota_K)$, on a 
$$
(\Phi_{\iota_{1,K}}(a_*(P)))(w)=(\Phi_{\iota_K}(P))(a^*(w)).
$$
Ici $$a_*:H^1(K,\ker\,\iota)\to H^1(K,\iota_{1,K})$$ est l'application induite par $a$ et 
$$a^*:F_K^*(\omega_{\ker\,\iota_{1,K}})\to F_K^*(\omega_{\ker\,\iota_{K}})$$ est l'application "image réciproque des formes différentielles par $a$". 
\label{compimp}
\end{comp}

\beginProof Soit $\lambda_P:P\to\Spec\,K$ un torseur sous $\ker\,\iota_K$. Comme plus haut, on sait (\cite[III.3.5.7]{Milne-Arith}) qu'il existe un unique morphisme 
$\sigma_P:\Spec\, K\to P$  tel que $\lambda_{P}\circ\sigma_P=F_K.$ Notons $\Gamma_1:=\ker\,\iota_1$. Soit $Q$ le quotient de $P\times_K\Gamma_{1,K}$ pour l'action de $\ker\,\iota_K=\Gamma_K$ donnée en coordonnées par 
la formule
\begin{equation}
(p,z)\mapsto(\gamma^{-1}(p),a(\gamma)+z)
\label{AcQ}
\end{equation}
On notera $$q:P\times_K\Gamma_{1,K}\to Q$$ l'application quotient (qui est finie et plate). L'action de $\Gamma_{1,K}$ par translation sur le deuxième facteur est compatible avec l'action décrite dans \refeq{AcQ} 
et $Q$ hérite ainsi d'une action de $\Gamma_{1,K}$ qui en fait un torseur sous $\Gamma_{1,K}$. Par définition, ce torseur représente $a(P)$. 
Soit $\lambda_Q:Q\to\Spec\,K$ le morphisme structural. On sait  qu'il existe un unique morphisme 
$\sigma_Q:\Spec\, K\to P$  tel que $\lambda_{Q}\circ\sigma_Q=F_K$ et on a encore une fois un diagramme commutatif

\begin{equation}
\xymatrix{
\Spec\,K\ar@/^3pc/[rrrr]^{=}\ar[r]\ar[ddr]_{\sigma_P\times\, 0} \ar@/_5pc/[dddrr]_{F_K}& F_K^*(\Gamma_K)\times_K F_K^*(\Gamma_{1,K})\ar[d]^{\sim}\ar[rr]^{\,\,\,\,\,\,(\gamma,\gamma_1)\mapsto a(\gamma)\cdot 
\gamma_1} & & F_K^*(\Gamma_{1,K})\ar[d]^{\sim}&\Spec\,K\ar[l]\ar[ddl]^{\sigma_Q}\ar@/^5pc/[dddll]^{F_K}\\
&F_K^*(P\times_K\Gamma_{1,K})\ar[rr]^{F_K^*(q)}\ar[d] & & F_K^*(Q)\ar[d]&\\
&P\times_K\Gamma_{1,K}\ar[rr]^{q}\ar[dr] & & Q\ar[dl]&\\
& &\Spec\,K & &
}
\label{ddone1}
\end{equation}
qui donne un diagramme commutatif
\begin{equation}
\xymatrix{
\sigma_Q^*(\Omega_{Q/K})\ar[rrr]\ar[d]^{\sim} & & & (\sigma_{P}\times 0)^*(\Omega_{P\times_K\Gamma_{1,K}/K})\ar[d]^{\sim}\\
F_K^*(\omega_{\Gamma_{1,K}}) \ar[rrr]^{a^*\oplus 0}\ar[d] & & & F_K^*(\omega_{\Gamma_K})\oplus F_K^*(\omega_{\Gamma_{1,K}}\ar[d])\\
\Omega_{F_K}\ar[rrr]^{=} & & & \Omega_{F_K}
}
\end{equation}
qui montre que 
l'application $F_K^*(\omega_{\Gamma_{1,K}})\to\Omega_{F_K}$ (qui n'est autre que 
$\Phi_{\iota_{1,K}}(P)$) 
est donnée par $\Phi_{\iota_K}(P)\circ a^*$, comme espéré.\endProof

\subsection{Le cas du morphisme de Frobenius relatif}

\label{ssFro}

{\it On suppose maintenant, et ce jusqu'à jusqu'à nouvel avis, que $$\iota=F_{\CA/V}:\CA\to \CA^{(p)}=\CB$$ est le morphisme de Frobenius relatif.} 

Voir par ex. \cite[3.2.4, p. 94]{Liu-Algebraic} pour la définition du morphisme de Frobenius relatif. Le schéma 
$\CA^{(p)}$ est par définition le changement de base de $\CA$ par $F_V$. 

Soit $x\in B(K)=A^{(p)}(K)$ et soit $P$ le torseur sous $\Gamma_K$ 
associé au $K$-schéma $\iota_K^*(x)$, où $x$ est vu comme un sous-schéma 
fermé réduit de $B$ et $\iota_K^*(x)$ est le changement de base de $x$  par $\iota_K$. Notons $u_x:P\hookrightarrow A$ l'immersion fermée naturelle. On a alors par construction une unique flèche 
$\tau_x:\Spec\, K\to P=\iota_K^*(x)$ telle que $\pi\circ u_x\circ \tau_x=F_K$. 
Si on note $K^{[1]}$ le corps $K$ vu comme une $k_0$-algèbre via l'application naturelle de 
$k_0$ dans $K$ précomposée par $F_{k_0}$, ceci induit une application
$$
\mu_x:\tau_x^*(u_x^*(\Omega_{A/k_0}))\to\Omega_{K^{[1]}/k_0}.
$$
Par ailleurs, on a une suite exacte sur $A$
$$
0\to\pi^*(\Omega_{K/k_0})\to \Omega_{A/k_0}\to\Omega_{A/K}\to 0
$$
et la composition de l'inclusion $\pi^*(\Omega_{K/k_0})\to \Omega_{A/k_0}$ 
avec $\mu_x$ est nulle car elle est induite par $F_{K}$. 
L'application $\mu_x$ induit donc une 
application 
$$
\tau_x^*(u_x^*(\Omega_{A/K}))\to \Omega_{K^{[1]}/k_0}
$$
et comme on a canoniquement
\begin{equation}
\tau_x^*(u_x^*(\Omega_{A/K})))\simeq \tau_x^*(\Omega_{\iota_K^*(x)/K})\simeq F_K^*(\omega_{\Gamma_K})
\label{alphacan}
\end{equation} 
et $$\Omega_{K^{[1]}/k_0}\simeq\Omega_{F_K},$$ 
on obtient une application
$$
\psi_x:F_K^*(\omega_{\Gamma_K})\to \Omega_{F_K}
$$
\begin{lemme}
On a $\psi_x=\Phi_{\iota_K}(P)$.
\label{lempsi}
\end{lemme}
\beginProof On a un diagramme

\smallskip
\hskip2cm
\xymatrix{
0\ar[r]&\tau_x^*(u^*_x(\pi^*(\Omega_{K/k_0})))\ar[dd]^{\simeq}\ar[r]&\tau_x^*(u^*_x(\Omega_{A/k_0}))\ar[dd]\ar[r]&
\tau_x^*(u^*_x(\Omega_{A/K}))\ar[ddl]\ar[r]\ar[d]^{\simeq}& 0\\
& & & \tau_x^*(\Omega_{\iota_K^*(x)/K})\ar[d]\ar[r]_{\simeq} & F^*_K(\omega_{\Gamma_K})\\
& F^*_K(\Omega_{K/k_0})\ar[r]^{=0}&\Omega_{K^{[1]}/k_0}\ar[r]^{\simeq}&\Omega_{F_K}\ar[r] &0}

\smallskip
dont les lignes sont exactes. On peut vérifier que ce diagramme commute 
en remplaçant $A$ par un voisinage affine de $x$ et en représentant 
les différentielles sous la forme $d(f)$, où $f$ est une fonction sur le voisinage affine. 
La commutativité du diagramme découle alors du fait que l'opération "image inverse d'une fonction" est fonctorielle et qu'elle commute avec l'opérateur $d(\cdot)$. Le lemme est une conséquence de la commutativité du diagramme.
\endProof

On note maintenant qu'on a aussi un isomorphisme naturel
\begin{equation}
F_K^*(\omega_A)\simeq\tau_x^*(u_x^*(\Omega_{A/K}))
\label{alphacan2}
\end{equation}
(puisqu'on a un isomorphismes canonique $\Omega_{A/K}\simeq\pi^*(\omega_A)$) qui donne en particulier un isomorphisme naturel
$$
\alpha_x:F_K^*(\omega_A)\stackrel{\sim}{\to} F_K^*(\omega_{\Gamma_K})
$$
via les isomorphismes \refeq{alphacan}.

{\it On notera le point important suivant.} Supposons le temps du prochain paragraphe que $$x\in\CB(V)\subseteq B(K).$$
Alors on peut remplacer $\Spec\, K$ par $V$ (resp. $A$ par $\CA$, resp. 
$B$ par $\CB$) dans les calculs précédents et on voit que 
$\alpha_x$ s'étend alors en un isomorphisme $$F_V^*(\omega_\CA)\stackrel{\sim}{\to} F_V^*(\omega_{\Gamma})$$ et $
\psi_x$ en une flèche $$F_V^*(\omega_{\Gamma})\to \Omega_{F_V}.
$$

\begin{lemme}
La flèche $\alpha_x$ ne dépend pas de $x\in A^{(p)}(K)$.
\label{alphalem} 
\end{lemme}
\beginProof 
On va démontrer cet énoncé par un argument de déformation utilisant la propreté de $A$. 
Soit $F=F_{A\times_K A\twp/A\twp}$ et soit $\Delta:A\twp\to A\twp\times_K A\twp$ le morphisme diagonal. 

On a un diagramme commutatif

\begin{equation}
\hskip4cm
\xymatrix{
A\times_K A\twp\ar[r]^{F} & A\twp\times_KA\twp\ar[d]_{\textrm sec. proj.}\\
A\twp\ar@{.>}[r]^<<<<{F_{A\twp}}\ar[u]^{\tau_\Delta}& A\twp\\
& A\ar[r]^<<<<<{F_{A/K}}\ar[luu] & A\twp\ar[d]\ar[luu]\\
&\Spec\ K\ar[r]^{F_K}\ar[u]^{\tau_x}\ar[uul]^{x}& \Spec\ K\ar@{.>}[luu]_>>>>>>>>>>>>>x
}
\label{I1}
\end{equation}
où les rectangles orthogonaux à la page sont cartésiens et 
où  $\tau_\Delta$ est défini de telle manière que 
$$
F\circ\tau_\Delta=\Delta\circ F_{A\twp}.
$$

On a aussi un diagramme commutatif

\begin{equation}
\hskip1.5
cm
\xymatrix{
A\times_K A\twp\ar[r]^{F}\ar[rd]\ar[ddr]_{\rm prem.\,proj.} & A\twp\times_KA\twp\ar[d]_{\textrm sec. pr.}\ar[ddr]^{\rm prem.\,proj.} \\
& A\twp\ar@{.>}[ddr]\\
& A\ar[r]^<<<<<{F_{A/K}}\ar[rd] & A\twp\ar[d]\\
&& \Spec\ K
}
\label{I2}
\end{equation}

où tous les carrés sont cartésiens. La flèche $A\to A\times_K A\twp$ 
du diagramme \refeq{I1} est une section de la flèche $A\times_K A\twp\to A$ du 
diagramme \refeq{I2}. De même la flèche $A\twp\to A\twp\times_K A\twp$ 
du diagramme \refeq{I1} est une section de la flèche $A\twp\times_K A\twp\to A\twp$ du 
diagramme \refeq{I2} et la flèche $x$ est une section de la flèche $A\twp\to\Spec\, K$.

Soit $$u_\Delta:F^*(\Delta_*(A\twp))\hookrightarrow A\times_K A\twp$$ 
l'immersion canonique dans $A\times_K A\twp$ du changement de base de la diagonale 
$\Delta_*(A\twp)$ par $F$. 

On a des isomorphismes canoniques
$$
\tau_\Delta^*(u_\Delta^*(\Omega_{A\times_K A\twp/A\twp}))
\simeq \tau_\Delta^*(\Omega_{F})\simeq 
F_{A^{(p)}}^*(\omega_{\ker\,F})
$$
et
$$
F_{A\twp}^*(\omega_{A\times_K A\twp/A\twp})\simeq\tau_\Delta^*(u_\Delta^*(\Omega_{A\times_K A\twp/A\twp}))
$$
analogues aux isomorphismes \refeq{alphacan} et \refeq{alphacan2}. 
On a donc un isomorphisme naturel
$$
\alpha:F_{A\twp}^*(\omega_{A\times_K A\twp/A\twp})\simeq F_{A^{(p)}}^*(\omega_{\ker\,F})
$$
analogue à l'isomorphisme $\alpha_x$. 

L'existence des diagrammes \refeq{I1} et \refeq{I2} et de leurs propriétés implique maintenant qu'on a un diagramme commutatif
\begin{equation}
\xymatrix{
F_K^*(\omega_A)\ar[r] & H^0(A\twp,F_{A\twp}^*(\omega_{A\times_K A\twp/A\twp}))\ar[r]\ar[d]^{H^0(\alpha)} & x^*(F_{A\twp}^*(\omega_{A\times_K A\twp/A\twp}))\ar[r]^<<<<<<{}\ar[d]^{x^*(\alpha)} & F_K^*(\omega_A)\ar[d]^{\alpha_x}\\
F_K^*(\omega_{\Gamma_K})\ar[r] & H^0(A\twp,F_{A^{(p)}}^*(\omega_{\ker\,F}))\ar[r] & x^*(F_{A^{(p)}}^*(\omega_{\ker\,F}))\ar[r]^<<<<<<<<<{}
& F_K^*(\omega_{\Gamma_K})
}
\label{alphaindep}
\end{equation}
où la composition des flèches dans la première ligne (resp. la deuxième ligne) 
est l'identité. Par ailleurs, toutes les flèches horizontales dans le diagramme 
\refeq{alphaindep} sont des isomorphismes car $A$ est propre, lisse et géométriquement connexe sur $\Spec\,K$. Le fait que $\alpha_x$ ne dépend pas de $x$ découle maintenant du fait que $H^0(\alpha)$ ne 
dépend pas de $x$. 
\endProof

\begin{cor}
La flèche $\alpha_x$ s'étend en un isomorphisme $F_V^*(\omega_\CA)\to F_V^*(\omega_{\Gamma})$.
\label{corai}
\end{cor}
\beginProof Puisque la flèche $\alpha_x$ ne dépend pas de $x$, nous pouvons 
sans restriction de généralité supposer que $x$ est la section nulle. 
Dans ce cas là, la flèche $\alpha_x$ s'étend en un isomorphisme $F_V^*(\omega_\CA)\to F_V^*(\omega_{\Gamma})$ par la remarque précédant le Lemme \ref{alphalem}.
\endProof

\begin{theor}
La flèche 
$$
\Phi_{\iota_K}(P):F_K^*(\omega_{\ker\ F_{A/K}})\to \Omega_{K^{[1]}/k_0}
$$
s'étend en une flèche 
$$
\phi_{\iota_K}(P):F_V^*(\omega_{\ker\ F_{\CA/C}})\to \Omega_{V^{[1]}/k_0}(\log\, E).
$$ 
\label{IItheor}
\end{theor}

Dans la preuve que l'on va lire, à la suite de L. Moret-Bailly on appellera "gros ouvert" un ouvert dont le complément est de codimension deux. On remplacera plusieurs fois de suite $V$ par l'un de ses gros ouverts pendant 
la démonstration. Ceci est licite car si $u:V_\circ\hookrightarrow V$ est un gros ouvert de $V$ et $F$ est un fibré cohérent localement libre sur $V$ alors on a $F\simeq u_*(u^*(F))$ puisque 
$V$ est normal. En particulier, si $F'$ est un autre fibré cohérent localement libre sur 
$V$ alors toute flèche $u^*(F)\to u^*(F')$ s'étend à une flèche $F\to F'$ de manière unique. 

\beginProof 
Nous rappelons d'abord quelques résultats démontrés dans \cite{FC-Degen}. 
On rappelle que $m$ est par hypothèse un entier $>2$.  
\begin{itemize}
\item[(1)] il existe un espace de modules $A_{g,m}$ pour les variétés abéliennes sur $k_0$ munies d'une structure de niveau $m$;
\item[(2)] il existe un schéma $A^*_{g,m}$ propre et lisse sur $k_0$ et une immersion ouverte $A_{g,m}\hookrightarrow A^*_{g,m}$, telle que le complément $D:=A^*_{g,m}\backslash A_{g,m}$ (vu comme sous-schéma fermé réduit) est un diviseur à croisements normaux dans $A^*_{g,m}$;
\item[(3)] il existe un schéma semiabélien $G$ sur $A^*_{g,m}$, qui étend le schéma abélien $f:Y\to A_{g,m}$ provenant de la propriété universelle de $A_{g,m}$;
\item[(4)] il existe un schéma régulier $\bar Y$ et un morphisme propre $\bar f:\bar Y\to A^*_{g,m}$ qui étend $f$ et le complément 
$F:=\bar Y\backslash Y$ (vu comme sous-schéma fermé reduit) est un diviseur à croisements normaux dans $\bar Y$; de plus
\item[(5)] il y a sur $\bar Y$ une suite exacte de faisceaux localement libres
$$
0\to \bar f^*(\Omega^1_{A_{g,m}/k_0}(\log D))\to\Omega^1_{Y/k_0}(\log F)\to\Omega^1_{Y/A_{g,m}}(\log F/D)\to 0,
$$
qui étend la suite habituelle de faisceaux localement libres
$$
0\to f^*(\Omega^1_{A_{g,m}/k_0})\to\Omega^1_{Y/k_0}\to\Omega^1_{Y/A_{g,m}}\to 0
$$
sur $A_{g,m}$. Enfin il y a un isomorphisme $$\Omega^1_{Y/A_{g,m}}(\log F/D)\simeq\bar f^*(\omega_G).$$ Ici $\omega_G:={\rm Lie}(G)^\vee$ est le fibré tangent (relativement à $A_{g,m}^*$) de $G$ restreint à $A_{g,m}^*$ via la section unité. 
\end{itemize}
Voir \cite[chap. VI, th. 1.1]{FC-Degen} pour la démonstration. 

La donnée de $A/K$ et de sa structure de niveau induit un morphisme $\phi:K\to A_{g,m}$ tel que $\phi^*Y\simeq A$, où l'isomorphisme préserve les structures de niveau. Appelons $\lambda:A\to Y$ 
le $k_0$-morphisme correspondant. Le critère valuatif de propreté implique que le morphisme 
$\phi$ s'étend en un morphisme $\bar\phi:V_\circ\to A^*_{g,m}$, où $V_\circ$ est un gros ouvert de $V$. 
On peut remplacer $V$ par $V_\circ$ (voir les remarques précédents la démonstration) et donc supposer que $\phi$ s'étend en un morphisme $\bar\phi:V\to A^*_{g,m}.$ 
Par l'unicité des modèles semiabéliens
(voir \cite[IX, Cor. 1.4, p. 130]{Raynaud-Faisceaux}), on a un isomorphisme naturel $\bar\phi^*(G)\simeq\CA$ (où $\bar\phi^*(G)$ est le changement de base de $G$ par $\bar\phi$) et on a donc une égalité ensembliste $\bar\phi^{-1}(D)=E$ 
et un isomorphisme canonique $\bar\phi^*(\omega_G)\simeq\omega_\CA$. 

De façon analogue, soit $\bar T_x:V_\circ\to\bar Y$ l'extension du morphisme $\lambda\circ\tau_x$ à un gros ouvert $V_\circ$ de $V$ obtenue via le critère valuatif de propreté. À nouveau on peut remplacer $V$ par $V_\circ$ et ainsi supposer que $\bar T_x$ est défini sur 
tout $V$. Par construction, on obtient maintenant une flèche 
$$
\bar T_x^*(\Omega^1_{\bar Y/k_0}(\log F))\to\Omega^1_{V^{[1]}/k}(\log E)
$$
et puisque la flèche induite 
$$
\bar T_x^*(f^*(\Omega^1_{A_{g,m}/k_0}(\log D)))=F^{*}_V\circ\bar\phi^*(\Omega^1_{A_{g,m}/k_0}(\log D))\to 
\Omega^1_{V^{[1]}/k_0}(\log E)
$$
s'annule (puisqu'elle s'annule génériquement), on obtient une flèche  
$$
\bar T_x^*(\Omega^1_{\bar Y/A^*_{g,m}}(\log F/D))=F^{*}_V\circ\bar\phi^*(\omega_G)=
F^{n,*}_V(\omega_\CA)\to \Omega^1_{V^{[1]}/k_0}(\log E),
$$
qui, en vertu du Lemme \ref{lempsi} et du Corollaire \ref{corai} est bien l'extension cherchée.\endProof

Le Théorème \ref{IItheor} montre déjà que le Théorème \ref{THprinc} est vérifié 
pour $P=\iota_K^*(x)$. 

\begin{rem}\label{goodhyp}\rm Soit  $L_0$ le corps de fonctions d'une variété quasi-projective 
sur $k_0$. Soit $C/L_0$ une variété abélienne. 
Alors il existe

- une extension de corps $L|L_0$ qui est finie et séparable;

- une variété projective $U$ sur $k_0$ dont le corps de fonctions est $L$;

- un schéma semiabélien $\CC$ sur $U$ tel que $\CC_L\simeq C_L$;

- un ouvert $U_\ab$ tel que $U\backslash U_\ab$ est un diviseur à croisements normaux 
et tel que $u\in U_\ab$ ssi $\CC_u$ est une variété abélienne.

Autrement dit, si on se donne une variété abélienne $C/L_0$ comme plus haut, il est possible après une extension finie et séparable $L$ de $L_0$ de trouver un modèle de $L$ qui est projectif et satisfait les hypothèses décrites dans le premier paragraphe de l'introduction.

Pour démontrer cette assertion, nous allons utiliser les notations 
de la démonstration du Théorème \ref{IItheor}, en particulier en ce qui concerne les espaces de modules de variétés abéliennes. 

Quitte à remplacer $L_0$ par une extension finie et séparable, on peut supposer 
qu'il existe $m\geq 3$ tel que $(p,m)=1$ et tel que $C[m](L_0)\simeq(\mZ/m\mZ)^{2\dim(C)}$. 
Soit $U_0$ une variété quasi-projective sur $k_0$ telle que $\kappa(U_0)=L_0$. 
Quitte à restreindre la taille de $U_0$, on peut supposer qu'il existe 
un $k_0$-morphisme $h:U_0\to A^*_{\dim(C),m}$ tel que 
$(h^*(Y))_{L_0}=C$. On choisit maintenant une immersion ouverte 
$U_0\hookrightarrow\bar{U}_0$, où $\bar U_0$ est une variété projective sur $k_0$. 
Considérons le graphe $\Gamma\subseteq U_0\times_{k_0}A_{\dim(C),m}^*$ de 
$h$. La première projection $\Gamma\to U_0$ est un isomorphisme; il existe donc un unique 
point $\eta\in \Gamma$ s'envoyant sur le point générique de $U_0$ par la première 
projection. On définit maintenant $U_1$ comme la clôture de Zariski de $\eta$ dans 
$\bar U_0\times_{k_0}A_{\dim(C),m}^*.$ La variété $U_1$ est par construction projective 
sur $k_0$ et la deuxième projection $U\to A_{\dim(C),m}^*$ la munit d'un schéma semiabélien $\CC_1$. Par ailleurs, elle est birationnelle à $U_0$ via la première projection 
$U_1\to\bar U_0$. 
\'Ecrivons $D_1$ pour l'image inverse de $D$ par la deuxième projection $U_1\to A_{\dim(C),m}^*.$  
On remplace maintenant $U_1$ par une altération $\alpha:U\to U_1$ génériquement étale  telle 
que $\alpha^{-1}(D_1)$ est un diviseur à croisements normaux. Une pareille altération 
existe par le fameux théorème \cite[Th. 4.1]{DJ-Alterations} de A.-J. de Jong. On rappelle 
qu'une altération est un $k_0$-morphisme propre et dominant de variétés sur $k_0$. On définit $L$ comme le corps de fonctions de $U$. 
La variété $U$ a toutes les propriété demandées.
\end{rem}

\subsection{Fin de la démonstration}

\label{ssFD}

{\it On relaxe maintenant la condition que $\iota=F_{\CA/V}$ et on suppose à nouveau seulement que $\iota$ est de hauteur un.}

La stratégie de la preuve est de se ramener au cas où $\iota$ est le morphisme de Frobenius relatif et où le torseur est donné par un point rationnel. Nous aurons besoins de quelques 
résultats préliminaires. 

\begin{lemme}
Pour que $P\in H^1(K,\Gamma_K)$ soit dans $\Sel_V(K,\iota)$ il faut et il suffit que pour tout point 
$v\in V$ de codimension un, il existe 

- un schéma intègre $U$ et un morphisme \'etale $U\to V$ dont l'image contient $v$;

- un point $u\in U$ tel que l'extension $\kappa(u)|\kappa(v)$ est de degré un;

- un élément $x\in A(L)$, où $L$ est le corps de fonctions de $U$, tel que $\delta_{\iota_L}(x)=P_L$ dans $H^1(L,\Gamma_L)$.

\label{lemA}
\end{lemme}
\beginProof Au vu de la suite exacte de Kummer \refeq{KSeq}, on voit que $P\in \Sel(K,\iota)$
ssi $t_v(P_{K_v})=0$ pour tout point de codimension un de $V$, en d'autres termes ssi pour tout point 
de codimension un $v$ de $V$ le torseur $T$ (à isomorphisme prêt) sous $A_{K}$ correspondant 
 à $t(P)$ a un point $K_v$-rationnel. Soit $v$ un point de codimension un de $V$. Par le théorème 
d'approximation de Greenberg, le torseur $T$ a un point $K_v$-rationnel ssi $T$ a un point $K_{v,h}$-rationnel, où $K_{v,h}$ est la henselisation de 
$K$ en $v$.  Le lemme découle maintenant de la définition de la henselisation.\endProof

Soit $L|K$ une extension finie et séparable de corps. En préparation du prochain lemme, on remarque que le 
diagramme 
\begin{center}
\hskip0.8cm
\xymatrix{
\Spec\, L\ar[r]^{F_L}\ar[d] & \Spec\, L\ar[d]\\
\Spec\, K\ar[r]^{F_K} & \Spec\, K
}
\end{center}
est cartésien. Pour vérifier cet énoncé, on note d'abord que le produit fibré 
$$\Spec\,L\times_{\Spec\, K}\Spec\, K$$ est topologiquement réduit à un point car il est fini et purement inséparable sur $\Spec\,L$ (car la propriété d'être fini et purement inséparable commute à tout changement de base). Par ailleurs, il est aussi lisse, puisqu'il est aussi le changement de base 
du morphisme lisse $\Spec\,L\to\Spec\, K$ à $\Spec\,K$. On conclut que  $\Spec\,L\times_{\Spec\, K}\Spec\, K$ est un corps et que le morphisme $\Spec\,L\times_{\Spec\, K}\Spec\, K\to\Spec\,L$ est 
fini, purement inséparable et de degré $\deg(F_K)$.  Par ailleurs, on a une flèche 
naturelle $\Spec\,L\to \Spec\,L\times_{\Spec\, K}\Spec\, K$ et comme 
les morphismes $F_L$ et $\Spec\,L\times_{\Spec\, K}\Spec\, K\to\Spec\,L$ sont tous deux de même degré, on conclut que cette flèche est un isomorphisme. 

\begin{lemme}
Soit $L|K$ une extension finie et séparable de corps. Soit 
$P$ un torseur sous $\Gamma_K$ et soit $P_L$ le torseur 
sous $\Gamma_L$ obtenu par changement 
de base. Soit 
$$
\beta_{L|K}:\Hom(F_K^*(\omega_{\Gamma_K}),\Omega_{F_K})\to 
\Hom(F_L^*(\omega_{\Gamma_L}),\Omega_{F_L})
$$
l'application obtenue en composant l'application de changement de 
base 
$$
\Hom(F_K^*(\omega_{\Gamma_K}),\Omega_{F_K})\to 
\Hom((F_K^*(\omega_{\Gamma_K}))_L,(\Omega_{F_K})_L)
$$
avec l'isomorphisme naturel $F_L^*(\omega_{\Gamma_L})\simeq (F_K^*(\omega_{\Gamma_K}))_L$ sur le premier facteur et l'isomorphisme naturel 
$(\Omega_{F_K})_L\to \Omega_{F_L}$ sur le deuxième facteur. 
Alors $$\Phi_{\iota_L}(P_L)=\beta_{L|K}(\Phi_{\iota_K}(P)).$$
\label{lemB}
\end{lemme}
\beginProof La vérification est élémentaire et nous l'omettons. \endProof

Soit $\beta:\CB\to \CB_1$ une isogénie de hauteur un telle que $\beta\circ\iota$ est aussi de hauteur un. On écrira $B_1:=\CB_{1,K}$. 
On a un diagramme commutatif

\smallskip
\begin{center}
\hskip1.5cm
\xymatrix{
 \ker\,\iota_K\ar[r]\ar[d] & A\ar[r]^{\iota_K}\ar[d]^{=} & B\ar[r]\ar[d]^{\beta_K}\ar[r]^{\delta_{\iota_K}\,\,\,\,\,\,\,\,\,\,\,\,\,\,} & H^1(K,\ker\,\iota_K)\ar[r]\ar[d] & H^1(K,A)\ar[d]^{=}\\
 \ker\,\beta_K\circ\iota_K\ar[r]\ar[d] & A\ar[r]^{\beta_K\circ\iota_K}\ar[d]^{\iota_K} & B_1\ar[r]\ar[d]^{=}\ar[r]^{\delta_{\beta_K\circ\iota_K\,\,\,\,\,\,\,\,}\,\,\,\,\,\,\,\,\,\,\,\,\,\,\,\,\,\,} & H^1(K,\ker\,\beta_K\circ\iota_K)\ar[r]\ar[d] & H^1(K,A)\ar[d]^{\iota_K^*}\\
  \ker\,\beta_K\ar[r] & B\ar[r]^{\beta_K} & B_1\ar[r]\ar[r]^{\delta_{\beta_K}\,\,\,\,\,\,\,\,\,\,\,\,\,\,} & H^1(K,\ker\,\beta_K)\ar[r] & H^1(K,B)
  }
 \end{center}
 \smallskip
 
L'application $H^1(K,\ker\,\iota_K)\to H^1(K,\ker\,\beta_K\circ\iota_K)$ est injective car son noyau est $H^0(K,\ker\,\beta_K)$, qui est nul puisque $\ker\,\beta_K$ est infinitésimal.  On a un 
diagramme tout semblable avec $K_v$ en place de $K$, pour toute valuation discrète $v$ sur $K$.
En examinant ces diagrammes et en tenant compte de la précédente remarque, on voit qu'on a une suite exacte
$$
0\to\Sel(K,\iota_K)\to\Sel(K,\beta_K\circ\iota_K)\to\Sel(K,\beta_K)
$$
Par ailleurs, les morphismes naturels de différentielles donnent lieu à un complexe
$$
0\to\omega_{\ker\,\beta}\to\omega_{\ker\,\beta\circ\iota}\to \omega_{\ker\,\iota}\to 0
$$
qui est exact en vertu de la classification des schémas en groupes de hauteur 
un par leurs algèbres de $p$-coLie (cf. \cite[Expos\'e VIIA, rem. 7.5]{SGA3-1}). Ceci suggère le 
 
\begin{lemme} On a un diagramme commutatif

\hskip 0.5cm
\xymatrix{
& 0\ar[d] & 0\ar[d]& 0\ar[d]\\
0\ar[r] & \Sel(K,\iota_K)\ar[r]\ar[d] & \Sel(K,\beta_K\circ\iota_K)\ar[r]\ar[d] & \Sel(K,\beta_K)\ar[d]\\
0\ar[r] & \Hom(F_K^*(\omega_{\ker\,\iota_K}),\Omega_{F_K})\ar[r] & 
\Hom(F_K^*(\omega_{\ker\,\beta_K\circ\iota_K}),\Omega_{F_K})\ar[r] & 
\Hom(F_K^*(\omega_{\ker\,\beta_K}),\Omega_{F_K})\ar[r] & 0
}

\smallskip
dont les colonnes et les lignes sont exactes.
\label{selkey}
\end{lemme}

\beginProof On doit seulement démontrer la commutativité du diagramme. Celle-ci est une 
conséquence immédiate du Complément \ref{compimp}.\endProof


Soit $x\in B(K)$ et soit comme avant $P$ le torseur sous $\Gamma_K$ 
associé au $K$-schéma $\iota_K^*(x)$, où $x$ est vu comme un sous-schéma 
fermé réduit de $B$ et $\iota_K^*(x)$ est le changement de base de $x$ par $\iota_K$.

\begin{theor} Soit $n\geq 3$. 
Supposons que $A[n](K)\simeq(\mZ/n\mZ)^{2g}$. 
La flèche 
$$
\Phi_{\iota_K}(P):F_K^*(\omega_{\Gamma_K})\to \Omega_{F_K}
$$
s'étend en une flèche 
$$
\phi_{\iota_K}:F_V^*(\omega_{\Gamma})\to\Omega_{V^{[1]}/k_0}(\log\, E)
$$
via l'isomorphisme canonique $\Omega_{F_K}\simeq\Omega_{K^{[1]}/k_0}$.
\label{THalth}
\end{theor}
\beginProof Soit $\beta:\CB\to \CB_1=\CA\twp$ l'unique isogénie de hauteur un telle  
que $$\beta_K\circ\iota_K=F_{A/K}.$$ On remarque maintenant qu'on a une suite exacte

\medskip
\centerline{$
0\to\Hom(F_V^*(\omega_{\ker\,\iota}),\Omega_{V^{[1]}/k_0}(\log\,E))\stackrel{a_V}{\to}\Hom(F_V^*(\omega_{\ker\,{\beta}\circ\iota}),\Omega_{V^{[1]}/k_0}(\log\,E))\stackrel{b_V}{\to} 
\Hom(F_V^*(\omega_{\ker\,{\beta}}),\Omega_{V^{[1]}/k_0}(\log\,E))\to 0
$}
\smallskip
qui se restreint à une suite exacte

\medskip
\centerline{
\xymatrix{
0\ar[r] & \Hom(F_K^*(\omega_{\ker\,\iota_K}),\Omega_{K^{[1]}/k_0})\ar[r]^{a_K} & 
\Hom(F_K^*(\omega_{\ker\,\beta_K\circ\iota_K}),\Omega_{K^{[1]}/k_0})\ar[r]^{b_K} & 
\Hom(F_K^*(\omega_{\ker\,\beta_K}),\Omega_{K^{[1]}/k_0})\ar[r] & 0}}

\medskip
 On sait par le Théorème 
\ref{IItheor} que $a_K(\Phi_{\iota_K}(P))$ se trouve dans
$$\Hom(F_V^*(\omega_{\ker\,{\beta}\circ\iota}),\Omega_{V^{[1]}/k_0}(\log\,E))\subseteq\Hom(F_K^*(\omega_{\ker\,\beta_K\circ\iota_K}),\Omega_{K^{[1]}/k_0}).$$ Soit 
$e\in \Hom(F_V^*(\omega_{\ker\,{\beta}\circ\iota}),\Omega_{V^{[1]}/k_0}(\log\,E))$ 
l'élément correspondant. Puisque $$b_K(a_K(\Phi_{\iota_K}(P)))=0$$ on sait que 
$b_V(e)=0$ et donc $e$ est l'image d'un élément $e'$ dans $\Hom(F_V^*(\omega_{\ker\,\iota}),\Omega_{V^{[1]}/k_0}(\log\,E))$. L'image de $e'$ dans 
$\Hom(F_K^*(\omega_{\ker\,\iota_K}),\Omega_{F_K})$ est par construction 
$\Phi_{\iota_K}(P)$, ce qui conclut la démonstration. 
\endProof

\beginProof (du Théorème \ref{THprinc}.) Nous aurons besoin du 

\begin{lemme}
Soit $S$ un schéma noethérien, intègre et normal et $H,J$ des faisceaux cohérents localement 
libres sur $S$. Soit $K$ le corps de fonctions de $S$. Soit $m:H_K\to J_K$. Supposons que pour tout point 
$s\in S$ de codimension un, il existe 

- un schéma intègre $U$ et un morphisme \'etale $U\to S$ dont l'image contient $s$;

- un point $u\in U$ tel que l'extension $\kappa(u)|\kappa(v)$ est de degré un;

- une extension de $m_L:H_L\to J_L$ à un morphisme $H_U\to J_U$, 
où $L$ est le corps de fonctions de $U$.

Alors il existe une extension de $m$ à un morphisme $H\to J$. 
\label{lemetext}
\end{lemme}
\beginProof (du Lemme \ref{lemetext}) On remarque tout d'abord qu'un morphisme 
$H\to J$ est un section globale d'un faisceau, à savoir $\underline{\Hom}_S(H,J)$. Par ailleurs, si 
une extension $m$ existe, elle est unique parce que $S$ est intègre; il en est de même pour une extension de $m$ à un ouvert de $S$. En mettant ensemble ces remarques, on voit qu'il suffit de montrer que pour tout point $s\in S$, on a 
\begin{equation}
m\in\Hom_{\CO_{S,s}}(H_{\CO_{S,s}},J_{\CO_{S,s}})\subseteq\Hom_K(H_K,J_K)
\label{extprop}
\end{equation}
Comme $S$ est normal, il suffit même de montrer que \refeq{extprop} est vérifié 
pour tout point $s\in S$ qui est de codimension un. Soit donc $s\in S$ un point de codimension un. Choisissons un base de $H_{\CO_{S,s}}$ et une base de $J_{\CO_{S,s}}$ (elles existent parce que $\CO_{S,s}$ est un anneau de valuation discrète). 
Une fois ces bases fixées, le morphisme $m$ est décrit par une matrice $M$ 
dont les coefficients sont des éléments de $K$. Soit maintenant $U\to S$ comme dans l'énoncé du lemme. Le morphisme $U\to S$ donne un diagramme 
commutatif

\smallskip
\begin{equation}
\xymatrix{
K \ar@{^{(}->}[r]\ar@/_3pc/@{_{(}->}[ddd] & \widehat{K}\ar@/^3pc/@{^{(}->}[ddd]\\
\CO_{S,s}\ar@{^{(}->}[r]\ar@{^{(}->}[d]\ar@{_{(}->}[u]& \widehat{\CO}_{S,s}\ar@{^{(}->}[d]\ar@{_{(}->}[u]\\
\CO_{U,u}\ar@{^{(}->}[r]\ar@{^{(}->}[d]& \widehat{\CO}_{U,u}\ar@{^{(}->}[d]\\
L \ar@{^{(}->}[r] & \widehat{L}
}
\label{lastdiag}
\end{equation}

Ici l'anneau $\widehat{\CO}_{S,s}$ (resp. $\widehat{\CO}_{U,u}$) est la complétion de l'anneau local $\CO_{S,s}$ (resp. ${\CO}_{U,u}$). Les corps 
$\widehat{K}$ et $\widehat{L}$ sont les anneaux de fractions de $\widehat{\CO}_{U,u}$ et $\widehat{\CO}_{U,u}$, 
respectivement. On notera que comme $\widehat{\CO}_{S,s}$ est un anneau de valuation 
discret, les carrés supérieurs et inférieurs du diagramme \refeq{lastdiag} sont cartésiens. 
Enfin, comme $U\to S$ est étale et que $\kappa(u)|\kappa(s)$ est un extension de degré un, 
on voit que la flèche $\widehat{\CO}_{S,s}\hookrightarrow \widehat{\CO}_{U,u}$ est un isomorphisme. 

Soit maintenant un coefficient $c\in K$ de la matrice $M$. 
Par hypothèse, l'image de $c$ dans ${L}$ est dans l'image de $\CO_{U,u}$. 
On déduit que l'image de $c$ dans $\widehat{K}$ est dans l'image de $\widehat{\CO}_{S,s}$ 
et (comme le carré supérieur est cartésien) on voit ainsi que $c$ est dans l'image de $\CO_{S,s}$, ce qu'on voulait montrer.\endProof

{\bf Fin de la démonstration du Théorème \ref{THprinc}.} On rappelle que la propriété d'un diviseur d'être à croisements normaux est par définition locale pour la topologie étale. 
Le Théorème \ref{THprinc} est donc une conséquence du Théorème \ref{THalth}, du 
Lemme \ref{lemB}, du Lemme \ref{lemA} et du Lemme \ref{lemetext}.\endProof

\end{document}